\documentclass[12pt]{amsart}
\usepackage{nccmath}

\usepackage{amsmath,amsthm,amssymb}
\usepackage{times}
\usepackage{enumerate}
\usepackage{verbatim}
\usepackage{xcolor}

\newtheorem{theorem}{Theorem}[section]
\newtheorem{corollary}[theorem]{Corollary}
\newtheorem{lemma}[theorem]{Lemma}

\newcommand{\be}{\begin{equation}}
\newcommand{\ee}{\end{equation}}
\newcommand{\ol}{\overline}

\newcommand{\goto}{\rightarrow}
\newcommand{\R}{\mathbb{R}}

\newcommand{\RNum}[1]{\uppercase\expandafter{\romannumeral #1\relax}}
\theoremstyle{definition}
\newtheorem{define}[theorem]{Definition}
\newtheorem{remark}[theorem]{Remark}

\newtheorem{claim}{Claim}
\usepackage[square,numbers]{natbib}
\date{\today \,(Last Typeset)}
\subjclass[2010]{Primary: 53C44, 53C21; Secondary:  53C42, 35J60.}
\keywords{Mean curvature flow, entire translating soliton, uniform 2-convexity,  bowl soliton,  convexity, fully nonlinear elliptic.}

\begin{document}
\setlength{\baselineskip}{1.2\baselineskip}

\title[ translating solitons]
{Convexity of 2-convex translating solitons to the mean curvature flow in $\R^{n+1}$}

\author{Joel Spruck}
\address{Department of Mathematics, Johns Hopkins University,
 Baltimore, MD 21218}
\email{js@math.jhu.edu}
\author{Liming Sun}
\address{Department of Mathematics, Johns Hopkins University,
 Baltimore, MD 21218}
\email{lsun@math.jhu.edu}

\begin{abstract} We prove that any complete immersed globally orientable uniformly 2-convex translating soliton $\Sigma \subset \R^{n+1}$ for the mean curvature flow is locally strictly convex. 
It follows that a uniformly 2-convex entire graphical translating soliton  in $\R^{n+1},\, n\geq 3 $ is the axisymmetric  ``bowl soliton''. 
\end{abstract}

\maketitle

\section{Introduction}
\label{sec0}
A solution of the mean curvature flow is a smooth one-parameter family
$\{\Sigma_t\}$ of hypersurfaces $\Sigma_t\subset \mathbb{R}^{n+1}$ with normal velocity equal to the mean curvature vector. A translating soliton for the mean curvature flow is the one  that  evolves purely by translation: $\Sigma_t = \Sigma +te$ for some fixed vector $e\in \mathbb{R}^{n+1}\backslash\{0\}$ and for all times $t$. In this case, the time slices are all congruent to $\Sigma$ and satisfy 
\begin{align}\label{eq:tran_eqn}
    \boldsymbol{H}=\langle \nu, e\rangle \nu=e^{\perp}
\end{align}
where $\nu$ is a choice of normal vector field for $\Sigma$ and $\boldsymbol{H}=-H\nu=-(\text{div}\, \nu)\nu$ is the corresponding mean curvature vector of $\Sigma$. We call $\Sigma$  a \textit{translator} in the direction $e$ for short. After a rotation, one can always assume $e=e_{n+1}$ after we normalize the speed to be one.

Translating solitons form a special class of eternal solutions for the mean curvature flow, that besides having their own intrinsic interest, are models of slow singularity formation. Therefore there has been a great deal of effort in trying to classify them in the most accessible case $H > 0$. 

For $n = 1$ the unique solution is the grim reaper curve $\Gamma: x_2 = \log \sec x_1$, $|x_1|<\pi$. For $n=2$, Wang \cite{Wang11} proved that any entire convex graphical translating soliton must be rotationally symmetric; this solution is  called the ``bowl soliton" \cite{Alt94}, \cite{CSS}. Moreover he showed that complete convex graphical translators that were not entire, necessarily live over strips. 
%The bowl soliton has the asymptotic expansion as an entire graph $u(x)=\frac{1}{2(n-1)}|x|^2-\log|x|+O(1)$. 
In \cite{SpruckLing18}, Spruck and Xiao proved that any complete immersed two-sided translating soliton in $\R^3$ with $H>0$ must be convex. Furthermore, they also classified the asymptotic behavior of possible solutions in a strip and conjectured the existence of a unique locally strictly convex translating soliton asymptotic to the ``tilted grim reaper'', $x_{3}=\lambda^2 \log{sec{\frac{x_1}{\lambda}}}+\sqrt{\lambda^2-1}\, x_2$, 
associated to any strip of width $\lambda \pi, \, \lambda>1$.

 Bourni et al.\,\cite{BLT18} proved the existence of such convex translators with the correct asymptotics (in fact in slabs of width greater than $\pi$ in any dimension $n\geq 2$). At about the same time,
Hoffman et al.\,\cite{Hoffman18} proved existence and uniqueness of locally strictly convex solitons in strips,  thus completing the classification of all mean convex translating solitons in $\mathbb{R}^3$, which consists of the standard grim reaper surface in a strip of width $\pi$, the tilted grim reaper in a strip of width $\pi \lambda>\pi$,
the locally strictly convex  ``delta wings''  asymptotic to the tilted grim reaper
at $\pm \infty$, in a strip of width $\pi \lambda>\pi$, and the  bowl soliton. 
 
In higher dimensions, Haslhofer \cite{Hasl15} proved the uniqueness of the bowl soliton in arbitrary dimensions under the assumption that the translating soliton $\Sigma$ is $\alpha$-noncollapsed and uniformly 2-convex. The $\alpha$-noncollapsed condition means that for each $P\in\Sigma$, there are closed balls $B^\pm$ disjoint from $\Sigma-P$ of radius at least $\frac{\alpha}{H(P)}$ with $B^+\cap B^-=\{P\}$. This condition figures prominently in the regularity theory for mean convex mean curvature flow \cite{HuiskenSinestrari99}, \cite{SW09}, \cite{White00}, \cite{White03}. The uniformly 2-convex condition (automatic if $n=2$) means that if $\kappa_1\leq \kappa_2\leq\cdots\leq \kappa_n$ are the ordered principle curvatures of $\Sigma$, then $\kappa_1+\kappa_2\geq \beta H$ for some uniform $\beta>0$. The $\alpha$-noncollapsed condition is a deep and powerful property of weak solutions of the mean convex mean curvature flow \cite{White03}, \cite{HK}, which implies that any complete $\alpha$-noncollapsed mean convex solition $\Sigma$
is convex with uniformly bounded second fundamental form. 
For some related results, see \cite{BrendleChoi17}, \cite{BrendleChoi18}, \cite{BL16}.
The purpose of  this paper is to extend  the work of \cite{SpruckLing18} to 2-convex translating solitons in all dimensions. The main result of this paper is the following
\begin{theorem}\label{thm1.1}
If $\Sigma^n$ is a complete immersed two-sided, uniformly 2-convex translator in $\mathbb{R}^{n+1},\, n\geq 3$, then $\Sigma$ is locally strictly convex.
\end{theorem}
Note that uniform 2-convexity plus convexity implies that there can be at most one zero curvature at a point of $\Sigma^n$.
Hence we are claiming in Theorem \ref{thm1.1}, that $\Sigma^n$ cannot
split off a line. For if $\Sigma^n=\Sigma^{n-1}\times \R$, then the uniform 2-convexity implies that the second fundamental form of $\Sigma^{n-1}$ satisfies $h_{ij}\geq \frac{\beta}n H$ and $\Sigma^{n-1}$ is a complete, locally strictly convex  graph over either a slab or all of $\R^{n-1}$. Hence by the main theorem of Hamilton \cite{Ham}, $\Sigma^{n-1}$ is compact, a contradiction.
As  a corollary, we obtain the following uniqueness theorem for the bowl soliton by appealing to Corollary 8.3 of a recent paper of Bourni, Langford and Tinaglia \cite{BLT19} that characterizes the bowl soliton.
\begin{corollary}\label{cor1.2} The bowl soliton is the unique uniformly 2-convex entire translator in $\mathbb{R}^{n+1},\, n\geq 3$.
\end{corollary}

For $n\geq 3$, the class of complete 2-convex translating solitons $\Sigma$ is rather restrictive. One elementary but important observation is that $|A^{\Sigma}|$ is uniformly bounded. In fact, all the principal curvatures have absolute value less than one by Lemma \ref{lem2.1}. In our proof of Theorem \ref{thm1.1},  we shall essentially utilize the property that $\kappa_1/(H-\kappa_1)$ satisfies a fully nonlinear elliptic equation. However, even though $\kappa_1$ is smooth in $\{\kappa_1<0\}$, the second fundamental form $h_{ij}=A(e_i, e_j)$ of $\Sigma$ may not be differentiable in a local curvature frame if there are positive curvatures with multiplicity.
To carry out our analysis, we utilize a special approximation $\mu^n(\kappa)$ of $\text{min}(\kappa_1,\ldots,\kappa_n)$ (section \ref{sec1}) which enjoys many good properties due to  our uniform 2-convex assumption. The $\mu^n(\kappa)$ depend on a  parameter $\delta$ and are defined recursively. Moreover, $\mu^n(\kappa)\to\kappa_1$ as $\delta\to 0$. We then apply the maximum principle  to  show that the infimum  of  $\mu^n/(H-\mu^n)$ cannot be achieved at a finite point, if $\delta$ is small enough. Thus the infimum must be achieved at infinity. As in \cite{SpruckLing18}, we apply the  Omari-Yau maximum principle  to a minimizing sequence of points $P_N$ tending to infinity. This argument is delicate and utilizes  the special properties of our approximation $\mu^n$. After  translating the $P_N$ back to the origin and passing to a subsequence, we again obtain a contradiction if $\delta$ small enough. Finally, this means $\mu^n$ can never be negative when $\delta$ is small enough. Therefore $\kappa_1\geq 0$ and $\Sigma$ is convex. 
%Recently, there are some work \cite{BrendleChoi17}, \cite{BrendleChoi18} which proved the uniqueness of ancient solution. And \cite{BLT19}......

The organization of the paper is as follows. In section \ref{sec1} we construct an approximation of the minimum function and derive several essential properties of it that will be needed in the proof of Theorem \ref{thm1.1}. In section \ref{sec2}, we extend the method of \cite{SpruckLing18} to prove the convexity of uniformly 2-convex translating solitons.

\section{Approximation of the minimum function}
\label{sec1}
In this section we refine the iterative approximation of $\text{min}\{x_1,\ldots, x_n\}$ of Heidusch \cite{Heidusch} and Aarons \cite{Aarons2005} so that it can be used in a  maximum principle argument to prove convexity of 2-convex translators for the mean curvature flow.

\begin{define}
\begin{enumerate}[(i)]
\item{The $\delta$-approximation to the function 
\[\min\{x_1,x_2\}\] is given by 
\[{\mu(x_1,x_2)=\mu^{2}\left(x_{1}, x_{2}\right)=\frac{x_{1}+x_{2}}{2}-\sqrt{\left(\frac{x_{1}-x_{2}}{2}\right)^{2}+\delta x_1x_2}} \]
for any  $\delta\in (0,\frac12)$.}

\item{  For $n\geq 3$, define inductively  the $\delta$-approximation to
\[\text{min}_n(x):=\min\{x_{1}, x_{2}, \cdots, x_{n}\},\text{ for }\, x\in \R^n\] by
\[\mu^n(x):=\frac{1}{n} \sum_{i=1}^{n} \mu(x_i,\mu^{n-1}(\bar{x}^i))\]
where $\bar{x}^i=(x_{1}, \cdots, x_{i-1}, x_{i+1}, \cdots, x_{n})\in \R^{n-1}$.}
\end{enumerate}
\end{define}

In order to accomondate our 2-convexity assumption and to show that the $\mu^n$ are well defined, we restrict $\mu^n$ to a convenient admissible domains $\mathcal{A}_n, \, \mathcal{A}^{-}_n,\, n\geq 2$, defined as follows:
\begin{align}
    \mathcal{A}_n=&\left\{x\in \mathbb{R}^{n}:   |x_j|\leq 1, \, \forall\,j ,\,\, x_k+x_l\geq \beta\sum_{i=1}^n x_i>0 \,\,\forall \, k\neq l \right\} ,\\
    \mathcal{A}^-_n=&\left\{x\in \mathcal{A}_n:\text{min}_n(x)\leq -\alpha \sum_{i=1}^nx_i\right\},
\end{align}
\begin{comment}
Change to |x_i|\leq 1 for all i
\end{comment}
for some fixed $\alpha, \beta\in (0,1)$. 
 It is easy to see that if $x\in \mathcal{A}_n^-$,   there is exactly one component  of $x$ with strictly negative minimum value   and all the other components are strictly positive. \\

\begin{lemma} \label{lem:nonempty}
For $n\geq 3$, a necessary and sufficient condition that  $\mathcal{A}_n^-\neq \emptyset$ is that
\begin{align}\label{eq2.10}
\beta \leq \frac{1-(n-2)\alpha}{n-1},\,\quad 0<\alpha<\frac1{n-2}.
\end{align}
\end{lemma} 
\begin{proof}
Without loss of generality, we may assume $x\in \mathcal{A}_n^-$ with $x_1\leq x_2\leq \cdots\leq x_n$. Then $x_1<0, \,\, x_2>0$ and 
\begin{align*}
    x_1+x_2\geq \beta\sum_{i=1}^nx_i>0\quad\text{and } \quad x_1\leq \frac{-\alpha}{1+\alpha}\sum_{i=2}^nx_i.
\end{align*}
or equivalently,
\[ x_1+x_2 \geq \frac{\beta}{1-\beta}\sum_{i=3}^nx_i \,, \,\,
x_1+\frac{\alpha}{1+\alpha}x_2 \leq \frac{-\alpha}{1+\alpha}\sum_{i=3}^nx_i,\]
which implies
\[\frac{1}{1+\alpha}x_2\geq \left[\frac{\beta}{1-\beta}+\frac{\alpha}{1+\alpha}\right]\sum_{i=3}^nx_i \geq (n-2)\left[\frac{\beta}{1-\beta}+\frac{\alpha}{1+\alpha}\right]x_2.\]
Since $x_2>0$, we arrive at 
\[\frac 1{1+\alpha}\geq  (n-2)\left[\frac{\beta}{1-\beta}+\frac{\alpha}{1+\alpha}\right] \]
which is equivalent to \eqref{eq2.10}. Conversely for any 
$0<\lambda\leq \frac{\alpha+1}{n-1}$, choose
\[\beta=\frac{1-(n-2)\alpha}{n-1},
\,x_1=-\frac{(n-1)\alpha \lambda}{1+\alpha},\,
x_2=x_3=\cdots=x_n=\lambda.\]
Then we have 
\[ m:=\sum_{i=1}^nx_i=\frac{(n-1)\lambda}{\alpha+1},\, x_1=-\alpha m,\, x_1+x_2=\beta m\]
so $x\in \mathcal{A}^{-}_n$.
\end{proof}

 The function $\mu$ has the following important properties we will need
 to analyze $\mu^n$.
%(see \cite{SAarons2005})

\begin{lemma}\label{lem2.0}
For any $0<\delta<\frac12$ and $(x_1,x_2)\in \mathcal{A}_2$,
\begin{enumerate}[{i.}]
\item{ $\,\mu$ is smooth and symmetric and if  $x_1,\,x_2>0$,  then $\mu(x)>0$. }
 \item{ $\,\mu$ is  monotonically increasing,  concave and satisfies
 \[ 0\leq \mu_{ x_i}\leq 1,  \,\, i=1, 2. \] }
 \item{ $\,\mu$ is homogeneous of degree $1$ and therefore satisfies $\sum_{i=1}^2 x_i\mu_{x_i}(x) =\mu(x)$.}
 \item{ $\,\mu(x)\leq \frac{1}{2}(x_1+x_2)$ and   $|\mu(x)-\min{(x_1,x_2)}|\leq 4\sqrt{\delta }\,(x_1+x_2)$.}
 \item{\, If $x\in \mathcal{A}_2^{-}$, assuming $x_1\leq x_2$, then
 \begin{align}\label{eq2.20}
 \frac12\left(1+\frac{1-2\delta+2\alpha(1-\delta)}{\sqrt{1+4\alpha(1+\alpha)(1-\delta)}}\right)\leq \mu_{x_1}<\frac12\left(1+\sqrt{1-\delta}\right)
 \end{align}
 and
 \begin{align}\label{eq2.25}
 \frac12\left(1-\frac{\alpha(1-2\delta)+1+\alpha}{\sqrt{1+4\alpha(1+\alpha)(1-\delta)}}\right)\leq \mu_{x_2}<\frac12 \left(1-\sqrt{1-\delta}\right).
 \end{align}}
 \end{enumerate}
 \end{lemma}
\begin{proof}
For i., it is easy to see that $\mu(x)$ is smooth except if $x_1=x_2=0$, and $\mu=0$ implies $\{x_1=0,x_2\geq 0\}$ or $\{x_1\geq 0,x_2=0\}$. If $x_1>0$ and $x_2>0$, then $\mu(x)>0$. 

For ii., a simple calculation show that
\begin{align}\label{eq2.30}
    \mu_{x_{1}}=\frac{1}{2}\left(1+\frac{(1-2\delta)x_{2}-x_1}{ \sqrt{\left(x_{1}-x_{2}\right)^{2}+4\delta x_1x_2}}\right) \in [0,1]
\end{align}
and $\mu_{x_1}=0$ if and only if $x_2=0$. Similarly
\begin{align}\label{eq2.35}
    \mu_{x_{2}}=\frac{1}{2}\left(1+\frac{(1-2\delta)x_{1}-x_{2}}{\sqrt{\left(x_{1}-x_{2}\right)^{2}+4\delta x_1x_2}}\right) \in [0,1]
\end{align}
and  $\mu_{x_2}=0$ if and only if $x_1=0$. Moreover,\begin{align}\label{eq2.36}
    (D^{2} \mu)=\frac{2\delta(1-\delta)}{\left[\left(x_{1}-x_{2}\right)^{2}+4\delta x_1x_2\right]^{\frac{3}{2}}} \left[ \begin{array}{cc}{-x_2^2} & {x_1x_2} \\ {x_1x_2} & {-x_1^2}\end{array}\right]
\end{align}
 is negative semi-definite. Hence  $\mu$ is increasing and concave. 

Statement  iii. is obvious as is the first part of statement iv.
 We verify the second statement of iv..
  By the homogeneity of $\mu$, it suffices  to restrict to  $x_1+x_2=1$ and
  we may also assume $x_1\leq x_2,\, -1\leq x_1\leq 1/2$. Then
   \begin{align*}
 |\mu-x_1| =&\frac12\left |\ \left((1-2x_1)-\sqrt{(1-\delta)(1-2x_1)^2+\delta}\, \right)\right |\\
 =&2\delta \left |\frac{x_1^2-x_1}{(1-2x_1)+\sqrt{(1-\delta)(1-2x_1)^2+\delta}}\right |\leq 4\sqrt{\delta},
  \end{align*}
  which gives $|\mu(x)-\min{(x_1,x_2)}|\leq 4\sqrt{\delta }$.
  
  To prove v., let $x_1=\text{min}_2(x)$ and observe that $x\in \mathcal{A}_2^{-}$ is equivalent to
  \[-1<t:=\frac{x_1}{x_2} \leq -\frac{\alpha}{1+\alpha}.\]
  Moreover a simple calculation shows that both
  \begin{align}\label{eq2.37}
  \mu_{x_1}=\frac12 \left(1+\frac{(1-2\delta)-t}{\sqrt{(t-1)^2+4\delta t}}\right),
  \end{align}
  and
  \begin{align}\label{eq2.38}
  \mu_{x_2}=\frac12\left(1+\frac{(1-2\delta)t-1}{\sqrt{(t-1)^2+4\delta t}}\right),
  \end{align}
are decreasing functions of $t$ on $(-1,-\frac{\alpha}{1+\alpha}]$ and 
then \eqref{eq2.20}, \eqref{eq2.25} follow by evaluation at the appropriate endpoints.
\end{proof}

We next show that $\mu^n$ is well-defined on $\mathcal{A}_n$ and has nice properties.
 \begin{lemma}\label{lem2.1} Let $x\in \mathcal{A}_n,\,n\geq 3$ and set $m:=\sum_{j=1}^n x_j$. Then
 \begin{enumerate}[i.]
 \item  $|\text{min}_{n}(x)|\leq 1-\frac{(n-2)\beta}{1-\beta} \,\,\text{if $\text{min}_{n}(x)<0$}$, 
 \item $\,  |\text{min}_n(x)| \leq \frac{1-\beta}{n-2}m,\,\, \max_i x_i \leq (1-\beta)m,$
 \item $x\in \mathcal{A}_n \Rightarrow \ol{x}^i\in \mathcal{A}_{n-1},\,\, i=1,\ldots,n,$
 \item $|\mu^n(x)-\text{min}_n(x)|\leq c(n)\sqrt{\delta}\,m,$
 \item $\text{For $\delta$ sufficiently small,}\,\, \ x\in \mathcal{A}_n \Rightarrow (x_i,\mu^{n-1}(\ol{x}^i))\in \mathcal{A}_2,\,\,i=1,\ldots,n,$
 \item For any $i$ such that $x_i\neq \text{min}_n(x)$, $ x\in \mathcal{A}^{-}_n \Rightarrow \ol{x}^i \in \mathcal{A}^{-}_{n-1},$
 \item $ \text{For $\delta$ sufficiently small,}\,\, x\in \mathcal{A}^{-}_n \Rightarrow (x_i,\mu^{n-1}(\ol{x}^i))\in \mathcal{A}^{-}_2, \,\,i=1,\ldots, n.
$
 \end{enumerate}
 \end{lemma}
 \begin{proof}
We prove i., ii. together.  Let $x\in \mathcal{A}_n\,$ and  assume
$\, x_1 \leq x_2 \leq \ldots \leq x_n$. Then
 \[\sum_{i=1}^n x_i=m\, \text{ with $x_1+x_2\geq \beta m$}~. \]
 In particular, 
 \[(1-\beta)(x_1+x_2)\geq \beta\sum_{i=3}^n x_i\geq (n-2)\beta x_2,\]
 so 
 \[|x_1|\leq \frac{1-(n-1)\beta}{1-\beta}x_2 \leq 1-\frac{(n-2)\beta}{1-\beta}\,\, \text{if $x_1<0$}.\] %add to statement of lemma
 %\textcolor{red}{this inequality is not right}
 In addition,  $\beta m+x_3+\ldots+x_n\leq m$. Then $|x_1|\leq x_2\leq x_3 \leq \frac{1-\beta}{n-2}m$ and it follows by induction that 
$x_n\leq (1-\beta)m$, proving i, ii.. Next observe that iii. is trivial if $i\geq 2$ or $i=1$ and $x_1\geq 0$. So assume $i=1$ and $x_1<0$. Then we must show
$x_2+x_3\geq \beta (\sum_{j=1}^n x_j-x_1)$ or equivalently
$\beta x_1+x_2+x_3 \geq \beta \sum_{j=1}^n x_j=\beta m$. Observe that
\begin{align*}
\beta x_1+x_2+x_3\geq (x_2-(1-\beta)x_1)&+x_1+x_2\geq \\
&(x_2-(1-\beta)x_1)+\beta \sum_{j=1}^n x_j>\beta m~,
\end{align*}
proving iii. (since we already know $|x_j|\leq 1$ for all $j$).  \\
We next prove iv. by induction; we  assume $x_1<x_2 \leq \ldots \leq x_n$. Since $\mu$ is homogenous of degree one, we will also assume $m=1$. The case $n=2$ is Lemma \ref{lem2.0} part iv. Now assume we have proved iv. for $n=2,\ldots,k-1$. Then by the monotonicity of $\mu$,   
\[
  \mu^k(x)=\frac1k\sum_{j=1}^k \mu(x_j, \mu^{k-1}(\ol{x}^j))
  =\frac1k(kx_1+O(\sqrt{\delta}))=x_1+O(\sqrt{\delta}),
 \]
% \textcolor{red}{I thought about this estimate for quite long time. It does not come easily to me. Basically we need to show $\mu(x_j,\mu^{k-1}(\bar x^j))=x_1+O(\sqrt{\delta})$ for any $i$. There are two many cases to consider. I am quite sure you acknowledge this or not. Anyway, we need to be careful about this point and do not have a circular proof.}\\
  which gives $|\mu^k(x)-\text{min}_k(x)|\leq c(k)\sqrt{\delta}$ completing the induction.\\
 We first  prove v. for $n=3$.
Let $x\in \mathcal{A}_3$ and note that if $x_1\geq 0$, then by part iii. and the definition of $\mu,\, 0\leq \mu(x_i, \ol{x}^i)<1,\, i=1,2,3$.  If $-1\leq x_1 <0$, then 
$$
\mu(x_1,x_2)-x_1 =\frac{x_2-x_1}2-\sqrt{\left(\frac{x_2-x_1}2\right)^2+\delta x_1x_2}>
 \frac{x_2-x_1}2-\frac{x_2-x_1}2=0.
 $$
 Hence
$\mu(x_1,x_2)\geq x_1\geq -1$ completing the proof for $n=3$.
Now suppose v. holds for $n=k-1$. Then
\[|\mu^{k-1}(x)|\leq \frac1{k-1}\sum_{j=1}^{k-1}|\mu(x_j, \mu^{k-2}(\ol{x}^j))|
\leq 1\]
since by induction, $|\mu^{k-2}(\ol{x}^j)|\leq 1$ and then $|\mu(x_j, \mu^{k-2}(\ol{x}^j)|\leq 1$ by the $n=3$ case. It remains to show
$x_i+\mu^{k-1}(\ol{x}^i) >0$ for all $i$. But again using part iv., we find (as in the proof of part iii.)
\[x_i+\mu^{k-1}(\ol{x}^i) \geq x_1+x_2-c(n) \sqrt{\delta}\, \sum_{j=1}^n x_j\geq (\beta-c(n)\sqrt{\delta}) \sum_{j=1}^n x_j >0\]
for $\delta$ small enough, completing the proof of part v.\\
 To prove vi.  it suffices by part iii. to show that for $i\geq 2,\,
 x_1\leq -\alpha(m-x_i)$ or equivalently,
 $x_1+\alpha m\leq \alpha x_i$.
 This holds trivially since the left hand side is nonpositive since $x\in \mathcal{A}^{-}_n$.\\
 Finally we prove vii. Again assume $m=1$. By part iv., 
 \begin{eqnarray*}
 &\text{min}(x_i, \mu^{n-1}(\ol{x}^i))\leq \text{min}_n(x)+c(n)\,\sqrt{\delta}=x_1+ c(n)\,\sqrt{\delta}\\
& x_i+\mu^{n-1}(\ol{x}^i)\leq x_n+x_1+c(n)\,\sqrt{\delta}.
 \end{eqnarray*}
% \textcolor{red}{should $c(n)$ be $c(n-1)$?}\\
Thus it suffices to show that 
 \[x_1+c(n)\sqrt{\delta}\leq -\alpha (x_n+ x_1+c(n)\,\sqrt{\delta})=
 -\alpha(1-(x_2+\ldots+x_{n-1})+c(n)\,\sqrt{\delta}).\]
 Since by hypothesis, $x\in \mathcal{A}^{-}_n,\, x_1\leq -\alpha m=-\alpha$,  we need only show that
 \[ x_2 \geq (1+\frac{1}{\alpha})c(n)\, \sqrt{\delta}.\]
  But  $x_1+x_2\geq \beta m=\beta$,
 hence $x_2 \geq\beta-x_1\geq  \beta+\alpha>(1+\alpha^{-1})c(n)\, \sqrt{\delta}\,$ for $\delta$ sufficiently small.
 \end{proof}

The following additional properties of the $\mu^n$ follow easily from
Lemmas \ref{lem2.0}, \ref{lem2.1} by induction.

\begin{corollary}\label{cor2.1}
For any $0<\delta<\frac12$, on $\mathcal{A}_n,\,n\geq 2$ we have:\\
\begin{enumerate}[{i.}]
\item{ $\mu^n$ is smooth and symmetric and if  $x_i>0$, $\forall\,i$ , then $\mu_n(x)>0$. }
 \item{$\mu^n$ is  monotonically increasing and concave and satisfies
 \[ 0\leq \mu^n_{ x_i}\leq 1\,\,\text{ for $\, i=1,\cdots,n$}. \] }
 \item{$\mu^n$ is homogeneous of degree $1$ and therefore $\sum_{i=1}^{n} x_i\mu^n_{x_i}(x) =\mu^n(x)$.}
 \item{$\mu^n(x)\leq \frac{1}{n}\sum_{j=1}^n x_j$.}
 \end{enumerate}
\end{corollary}

\begin{comment} 
This has been moved to Lemma \ref{lem2.0} part (5)
\begin{remark}\label{rem2.1}
If $x=(x_1,x_2)\in \mathcal{A}_2^{-}$, then $x_1=\text{min}_2(x)$ satisfies
\[ -x_2<x_1\leq -\frac{\alpha}{1+\alpha} x_2,\,\, \text{some $\alpha >0$},\]
and we can improve the estimates for 
$\mu_{x_1}$ and $\mu_{x_2}$. For then $\frac{x_1}{x_2}\in (-1, -\frac{\alpha}{1+\alpha})$ and $\mu_{x_1}$ is decreasing in $x_1/x_2$ in this range and it follows from \eqref{eq2.35} that
\begin{align}
    \frac{1}{2}+\frac{\frac{-\alpha}{\alpha+1}-1+2\delta}{4\sqrt{\left(\frac{-\frac{\alpha}{1+\alpha}-1}{2}\right)^2-\delta\frac{\alpha}{1+\alpha}}}\leq\mu_{x_1}\leq\frac12(1+\sqrt{1-\delta}).
\end{align}
Similarly $\mu_{x_2}$ is also decreasing in $x_1/x_2$ in this range and\eqref{eq2.35} shows
\begin{align}\label{eq2.40}
    \frac12+\frac{(1-2\delta)\frac{-\alpha}{\alpha+1}-1}{4\sqrt{\left(\frac{-\frac{\alpha}{1+\alpha}-1}{2}\right)^2-\delta\frac{\alpha}{1+\alpha}}}\leq \mu_{x_2}\leq\frac12(1-\sqrt{1-\delta}) 
\end{align}
\end{remark}
\end{comment}

Next we shall prove some important properties related to the derivatives of $\mu^n$.
\begin{lemma}\label{lem2.2}
For $n\geq 2$, there exist $\delta_n=\delta_n(\alpha,\beta,n)$ such that the following properties hold for $0<\delta<\delta_n$:
\begin{enumerate}[{i.}]
    \item If $x\in \mathcal{A}_n$ and $x_i<x_j$, then $\mu^n_{x_i}>\mu^n_{x_j}$.
    \item If $x\in \mathcal{A}_n^-$, then
    \begin{align*}
        \mu^n_{x_i}(x)\to\begin{cases}
        1&\text{if }x_i=\min_n(x)\\
        0&\text{otherwise }
        \end{cases}
    \end{align*} uniformly as $\delta\to 0^+$.
    \item If $x\in \mathcal{A}_n^-$, then $\mu^n_{x_i x_i}(x)<0$.
\end{enumerate}
\end{lemma}
\begin{proof}
We first verify properties i., ii., iii. for $n=2$. From the explicit formulas \eqref{eq2.30} and \eqref{eq2.35}, property i is obvious. If $(x_1,x_2)\in \mathcal{A}_2^-$, then as in Lemma \ref{lem2.0} part v for $x_1=\text{min}_2(x),\, t:=\frac{x_1}{x_2}\in (-1, -\frac{\alpha}{1+\alpha}]$ and one sees that
$\mu_{x_1} \goto1$ by \eqref{eq2.37} and $\mu_{x_2}\goto 0$ by \eqref{eq2.38} as $\delta\to 0^+$. Therefore property ii. is established and property iii. follows from \eqref{eq2.36}. 

For $n\geq 3$ we prove properties i., ii., iii. by induction. Assume they hold for  $\mu^{n-1}$. 
For i., it is enough to show that if $x_1<x_2$, then $\mu^n_{x_1} >\mu^n_{x_2}$. Write
\begin{align}\label{eq2.50}
    \mu^n(x)=d^n(x)+r^n(x)
\end{align}
where 
\begin{align*}
    d^n(x)=&\frac1n \mu(x_1,\mu^{n-1}(\bar{x}^1))+\frac1n \mu(x_2,\mu^{n-1}(\bar{x}^2))\\
  r^n(x)=&\frac{1}{n}\sum_{i=3}^n \mu(x_i, \mu^{n-1}(\bar{x}^i)).
\end{align*}
One sees easily that $d^n(x)$ and $r^n(x)$ are both symmetric and concave in  the variables $x_1$ and $x_2$. Consequently $d^n_{x_1}\geq d^n_{x_2}$ and $r^n_{x_1}\geq r^n_{x_2}$. 
%(for example, see \cite{Spruck}). 
Furthermore,  
\begin{align} \label{eq2.60}
    &r^n_{x_1}(x)-r^n_{x_2}(x)=\frac{1}{n}\sum_{i=3}^n \mu_{y_2}(x_i,\mu^{n-1}(\bar{x}^i))[\mu^{n-1}_{x_1}-\mu^{n-1}_{x_2}] (\bar{x}^i).
\end{align}
Here $\mu_{y_2}$ means the partial derivative with respect to the second variable.
It follows from Lemma \ref{lem2.1} that $\bar{x}^i\in \mathcal{A}_{n-1}$ and  $(x_i,\mu^{n-1}(\bar{x}^i))\in \mathcal{A}_2$.
 Therefore, 
 $$
 \mu_{y_2}(x_i,\mu^{n-1}(\bar{x}^i))\geq 0,
 $$
 with equality for such  an $i$ only when $x_i=0$.  By our assumption on $\mu^{n-1}$, one knows 
 $$
 \mu^{n-1}_{x_1}(\bar{x}^i)>\mu^{n-1}_{x_2} (\bar{x}^i).
 $$
  Since $x\in\mathcal{A}_n\,$,  $x$ cannot have two entries which are  both equal to $0$. Therefore $r^n_{x_1}(x)-r^n_{x_2}(x)>0$ except possibly if $n=3$ and $x_3=0$. However, in this case one must have $x_1, x_2>0$ and consequently $\mu^3(x)=\mu(x_1,x_2)$, a case already verified. The proof of property i. is complete.

To prove property ii., we may assume $x_1=\text{min}_n(x)$. Since $x\in \mathcal{A}_n^{-},\,\,x_1<0$. Note that
\begin{align*}
    \mu^n_{x_1}(x)=\frac{1}{n}\mu_{y_1}(x_1,\mu^{n-1}(\bar{x}^1))+\frac{1}{n}\sum_{i=2}^n\mu_{y_2}(x_i,\mu^{n-1}(\bar{x}^i))\mu_{x_1}^{n-1}(\bar{x}^i).
\end{align*}
 It follows from Lemma \ref{lem2.1} that $(x_i,\mu^{n-1}(\bar{x}^i))\in \mathcal{A}_2^-$ for any $1\leq i\leq n$ and $\bar{x}^i\in \mathcal{A}_{n-1}^-$ for any $i\geq 2$. Again by our assumption on $\mu^{n-1}$, we see that $\mu^n_{x_1}\to 1$ as $\delta\to 0$. One can prove similarly that $\mu^n_{x_i}\to 0$ for $i\geq 2$.

To prove iii., we split $\mu^n$ as 
\[\mu^n(x)=\frac{1}{n}\mu(x_i,\mu^{n-1}(\bar{x}^i))+\frac{1}{n}\sum_{j\neq i}\mu(x_j,\mu^{n-1}(\bar{x}^j))\]
The last term is concave in $x_i$ and its second pure derivative in $x_i$ is non-positive. Thus we need to show
\begin{align}
    \mu_{x_i x_i}(x_i,\mu^{n-1}(\bar{x}^i))=\mu_{y_1y_1}(x_i,\mu^{n-1}(\bar{x}^i))<0.
\end{align}
However, this can be seen from $(x_i,\mu^{n-1}(\bar{x}^i))\in\mathcal{A}_2^-$. The proof by induction is now complete. 
\end{proof}
Recall from Lemma \ref{lem2.0} part v. that if $x\in \mathcal{A}_2^{-}$ with $x_1<0<x_2$,
 $$
\mu_{x_2}\geq \Lambda=\Lambda(\alpha,\delta):= \frac12\left(1-\frac{\alpha(1-2\delta)+1+\alpha}{\sqrt{1+4\alpha(1+\alpha)(1-\delta)}}\right)  >0.
$$
\begin{lemma}\label{lem2.3}
Suppose $n\geq 3$. There exists $\delta_n=\delta_n(\alpha,\beta,n)$ such that if $\,0<\delta<\delta_n$ and $x\in \mathcal{A}_n^-$, then 
    \begin{align}
        \frac{\mu^n_{x_i}-\mu^n_{x_j}}{x_i-x_j}\leq -\frac{\Lambda^{n-2}}{2\cdot n!}\,\frac{1}{|x_i-x_j|+\sqrt{x_ix_j}}
    \end{align}
  for any $x_i\neq x_j$ and $x_i>0,x_j>0$.
\end{lemma}
\begin{proof}
We first examine the case $n=2$. If $x_1\neq x_2$ and $x_1,x_2>0$, then it follows from \eqref{eq2.30} and \eqref{eq2.35} that
\begin{align}
    \frac{\mu_{x_1}-\mu_{x_2}}{x_1-x_2}=-\frac{1-\delta}{2 \sqrt{\left(x_{1}-x_{2}\right)^{2}+4\delta x_1x_2}}\leq -\frac{1}{4|x_1-x_2|+4\sqrt{x_1x_2}}
\end{align}
provided $\delta<\frac14$. Now consider $n\geq 3$, using the decomposition $\mu^n=d^n+r^n$ as in \eqref{eq2.50}. Then 
\[\frac{\mu^n_{x_i}-\mu^n_{x_j}}{x_i-x_j}\leq \frac{r^n_{x_i}-r^n_{x_j}}{x_i-x_j} \]
We may assume $x_1\leq x_2\leq \cdots\leq x_n$ and $x_1<0$ since $x\in \mathcal{A}_n^-$. Notice we are assuming $x_i,x_j>0$, which implies $i>1$, $j>1$. 

Supposing $n=3$ and recalling \eqref{eq2.60}, one obtains
\begin{align*}
    \frac{r^3_{x_2}-r^3_{x_3}}{x_2-x_3}\leq& \frac{1}{3}\mu_{y_2}(x_1,\mu(\bar{x}^1))\cdot \frac{\mu_{x_2}(\bar x^1)-\mu_{x_3}(\bar x^1)}{x_2-x_3}\\
    \leq& -\frac{\Lambda}{2\cdot 3!}\frac{1}{|x_2-x_3|+\sqrt{x_2x_3}}
\end{align*}
We use induction to prove the cases $n>3$. Let $l$ be an number in $\{2,\cdots,n\}$  with $l\neq i, j$. Again using \eqref{eq2.60},
\begin{align*}
    \frac{r^n_{x_i}-r^n_{x_j}}{x_i-x_j}\leq \frac{1}{n}\mu_{y_2}(x_l,\mu^{n-1}(\bar{x}^l))\cdot \frac{\mu^{n-1}_{x_i}(\bar x^l)-\mu^{n-1}_{x_j}(\bar x^l)}{x_i-x_j}
\end{align*}
By Lemma \ref{lem2.2}, $(x_l,\mu^{n-1}(\bar{x}^l))\in \mathcal{A}_2^-$, so it follows from \eqref{eq2.35} that 
\begin{align*}
    \frac{r^n_{x_i}-r^n_{x_j}}{x_i-x_j}\leq \frac{\Lambda}{n}\, \frac{\mu^{n-1}_{x_i}(\bar x^l)-\mu^{n-1}_{x_l}(\bar x^1)}{x_i-x_j}\leq -\frac{\Lambda^{n-2}}{2\cdot n!}{\frac{1}{|x_i-x_j|+\sqrt{x_ix_j}}},
\end{align*}
where the last inequality follows from our induction  hypothesis and $\bar x^l \in \mathcal{A}_{n-1}^-$.
                                                                                                                                                                                                                                                                                                                                                                                                                                                                             
\end{proof}
\section{2-convex translators}\label{sec2}
Suppose $e_{n+1}$ is the direction of the translation. The mean curvature $H= \text{tr}A$ and second fundamental form $A=(h_{ij})$ satisfy the following equations, for instance see \cite{MSS15}
\begin{align}\label{eq:H-eqn}
    \Delta H+\nabla_{e_{n+1}}H+|A|^2H=0,\\
    \Delta A+\nabla_{e_{n+1}} A+|A|^2A=0.\label{eq:A-eqn}
\end{align}
Define $L=\Delta +\nabla_{e_{n+1}}$, which is so called drift laplacian. Suppose $\kappa_1,\cdots,\kappa_n$ are  the principle curvatures of $\Sigma$ and $\tau_1,\cdots, \tau_n$ is a smooth orthonormal frame. We write $\mu^n=\mu^n(\kappa_1,\cdots,\kappa_n)$ for short.

Suppose $\Sigma$ is uniformly 2-convex. More precisely, assuming $\kappa_1\leq \kappa_2\leq\cdots\leq \kappa_n$, we have $\kappa_1+\kappa_2\geq \beta H>0$ for some $\beta>0$. Then $(\kappa_1,\cdots, \kappa_n)\in \mathcal{A}_n$. By Corollary  \ref{cor2.1},  $\mu^n$ is smooth on $\Sigma$.
Since $\mu^n$ is  a symmetric function of  the principle curvatures, we can express $\mu^n$ as a fully nonlinear equation of the second fundamental form, namely $\mu^n(\kappa_1,\cdots,\kappa_n)=F((h_{j}^i))$. Here we point out that since $\mu^n$ is a smooth function and symmetric on its arguments, then $F$ is smooth on the second fundamental forms. Denote $F_{ij}=\frac{\partial F}{\partial h_{ij}}$ and $F_{ij,rs}=\frac{\partial^2F}{\partial h_{ij}\partial h_{rs}}$. Then a standard calculation gives
\begin{align*}
    \Delta \mu^n=& F^{ij}h_{ijkk}+F^{ij,rs}h_{ijk}h_{rsk}\\
    =&F^{ij}(-|A|^2h_{ij}-\nabla_{e_{n+1}} h_{ij})+F^{ij,rs}h_{ijk}h_{rsk}\\
    =&-|A|^2F^{ij}h_{ij}-\nabla_{e_{n+1}} \mu +F^{ij,rs}h_{ijk}h_{rsk}
\end{align*}
where we have used \eqref{eq:A-eqn}. Recall that by the definition of $L$, 
\begin{align*}
    L\mu^n=-|A|^2F^{ij}h_{ij}+F^{ij,rs}h_{ijk}h_{rsk}.
\end{align*}
Define $Q_\delta=\frac{\mu^n}{H-\mu^n}$. Then a simple calculation gives
\begin{align*}
    \Delta Q_\delta=\frac{H\Delta\mu^n-\mu^n\Delta H}{(H-\mu^n)^2}-2\left\langle\frac{\nabla (H-\mu^n)}{H-\mu^n}, \nabla Q_\delta\right\rangle
\end{align*}
which means
\begin{align}\label{eq:LQ}
    LQ_\delta+2\frac{\left\langle{\nabla (H-\mu^n)}, \nabla Q_\delta\right\rangle}{H-\mu^n}=\frac{HL\mu^n-\mu^n LH}{(H-\mu^n)^2}.
\end{align}
The previous calculation of $L\mu^n$ and \eqref{eq:H-eqn} shows that
\begin{align}\label{eq:LQ2}
    HL\mu^n-\mu^n LH={H|A|^2 (\mu^n -F^{ij}h_{ij})+HF^{ij,pq}h_{ijk}h_{pqk}}
    =HF^{ij,pq}h_{ijk}h_{pqk},
\end{align}
since $F^{ij}h_{ij}= \mu^n$ by the homogeneity property of $\mu^n$. 
Furthermore, the concavity of $\mu^n$ implies $F^{ij,rs}$ is negative definite. Since $H>0$, it follows from \eqref{eq:LQ}, \eqref{eq:LQ2} that 
\begin{align}\label{eq:Q-ellip}
    LQ_\delta+2\frac{\left\langle{\nabla (H-\mu^n)}, \nabla Q_\delta\right\rangle}{H-\mu^n}=\frac{HF^{ij,pq}h_{ijk}h_{pqk}}{(H-\mu^n)^2}\leq 0
\end{align}

We can restate our main result as 
\begin{theorem}\label{thm3.1}
If  $\, \Sigma^n$ is a uniformly 2-convex translator with $n\geq 3$, then 
\[\lim_{\delta\to 0^+}\inf_{\Sigma}Q_\delta\geq 0.\]
As a consequence, $\Sigma$ must be convex.
\end{theorem}
\begin{proof}
We will prove the theorem by contradiction. Recall that $\mu^n<\frac1nH$ and $H>0$. It is easy to see that $Q_\delta> -1$ on $\Sigma$. Assume $\inf_{\Sigma}Q_\delta<-\varepsilon_0<0$ for any $\delta>0$ small. From now on, we choose $\tau_1,\tau_2,\cdots,\tau_n$ to be principle directions corresponding to the ordered principle curvatures $\kappa_1\leq\cdots\leq \kappa_n$.

Suppose $Q_\delta$ attains its infimum at some interior point $P\in \Sigma$. Applying the strong maximum principle to \eqref{eq:Q-ellip} yields $Q_\delta\equiv Q_\delta(P)<-\varepsilon_0<0$ and $F^{ij,pq}h_{ijk}h_{pqk}=0$. In particular, $\mu^n=\frac{Q_\delta}{1+Q_\delta}H$. It follows from Lemma \ref{lem2.1} (iv)  that,
\begin{align}\label{eq:k1aH}
    \kappa_1\leq \mu^n+c({n})\sqrt{\delta} H= \frac{Q_\delta}{1+Q_\delta}H+c({n})\sqrt{\delta}H\leq -\alpha H
\end{align}
for some $\alpha>0$ small, if $\delta$ is taken small enough. We fix $\alpha,\beta>0$ small  enough  that \eqref{eq2.10} holds. Consequently, $(\kappa_1,\cdots,\kappa_n)\in \mathcal{A}_{n-1}$ for such $\alpha$ and $\beta$ and all points on $\Sigma$. 

\begin{claim} \label{cl1} If $Q_\delta$ has an interior infimum, then $h_{ijk}\equiv 0$, that is $\nabla A\equiv0$.
\end{claim}
\begin{proof}
With the notation that $\mu^n_i=\frac{\partial\mu^n }{\partial \kappa_i}$ and $\mu^n_{ij}=\frac{\partial^2 \mu^n}{\partial\kappa_i\partial\kappa_j}$,  it is well known that $F^{ij,pq}h_{ijk}h_{pqk}$ can be calculated (see for example, \cite{Spruck} or \cite{BenA}):
\begin{align}\label{eq:second-deri}
    F^{ij,pq}h_{ijk}h_{pqk}=\sum_{i,j,k}\mu^n_{ij}h_{iik}h_{jjk}+\sum_{\{i,j:\kappa_i\neq \kappa_j\}}\sum_{k}\frac{\mu^n_i- \mu^n_j}{\kappa_i-\kappa_j}h_{ijk}h_{ijk}.
\end{align}
By concavity, both  terms on the right hand side  of \eqref{eq:second-deri}
are nonpositive. Since $Q_\delta$ is constant, then \eqref{eq:Q-ellip} implies  $F^{ij,pq}h_{ijk}h_{pqk}=0$.  Therefore both of the above terms must be 0. Because of Lemma \ref{lem2.2}, 
\[\frac{\mu_i^n-\mu^n_j}{\kappa_i-\kappa_j}< 0\quad \text{for }\kappa_i\neq \kappa_j.\]
Thus  necessarily, for each $i$ and  $j$ such that $\kappa_i\neq \kappa_j$
\begin{align*}
    h_{ijk}=0,\quad \forall\, k.
\end{align*}

However, if $i\neq j$ and $\kappa_i=\kappa_j$ , then $h_{ijk}=0$ for any $k$. Indeed, for a dense open set of points, we can choose a smooth  orthonormal frame $\tau_1, \ldots, \tau_n$ of eigenvectors for the ordered principal curvatures $\kappa_1\leq ...\leq \kappa_n$, see Theorems 2 and 3 of \cite{Singley1975}. Then (see formula (9) of \cite{BLT18}),
\begin{align}\label{eq:hijk=0}
\begin{split}
    0&=\tau_k h_{ij}=\nabla_{\tau_k} h_{ij}+h(\nabla_{\tau_k}\tau_i,\tau_j)+h(\tau_i,\nabla_{\tau_k}\tau_j)\\
    &=h_{ijk}+\Gamma_{ki}^j \kappa_j+\Gamma_{kj}^i \kappa_i
  =h_{ijk}+\Gamma_{ki}^j(\kappa_j-\kappa_i)=h_{ijk}.
\end{split}
\end{align}
Summing up the above analysis, one has  $h_{ijk}\equiv 0$ for $i\neq j$ and any $k$, for a dense open set of points of $\Sigma$. Since $A=(h_{ij})$ is a Codazzi tensor, it remains to show $h_{iii}=0$ for all $i$. Recall that we also know $\sum_{i,j,k}\mu^n_{ij}h_{iik}h_{jjk}=0$, which  now  reduces to $\sum_{i}\mu^n_{ii}h_{iii}^2=0$. Since $\mu_{ii}^n<0$  as a consequence of  Lemma \ref{lem2.2} iii., it follows that $h_{iii}=0$ for all $i$.  Claim \ref{cl1} is established for a dense open set of points and hence all points by continuity.
\end{proof}

Now if $\nabla A\equiv 0$, then $\nabla H\equiv 0$. It follows that
\[0=\nabla_{\tau_l}H=\kappa_l\langle \tau_l,e_{n+1}\rangle. \]
Since $\kappa_l\neq 0$ for any $1\leq l\leq n$, then the unit normal of $\Sigma$ is $\nu=e_{n+1}$ and $H\equiv 1$ which is impossible. Therefore $Q_\delta$ cannot have an interior infimum. 

Therefore the infimum of $Q_\delta$ is achieved at infinity and we can apply the Omori-Yau maximum principle. That is,  there exists a sequence$P_{\delta,N}\to \infty$ such that 
\begin{align}\label{eq:OY}
    Q_\delta(P_{\delta,N})\to \inf_{\Sigma}Q_\delta,\quad |\nabla Q_\delta(P_{\delta,N})|<\frac{1}{N},\quad \Delta Q_\delta(P_{\delta,N})>-\frac{1}{N}
\end{align}
Moreover, we can perturb the $P_{\delta,N}\to \infty$ slightly so that \eqref{eq:second-deri} and 
\eqref{eq:hijk=0} are also satisfied.

%\textcolor{red}{Since $\Sigma$ is uniformly 2-convex, one can choose $\alpha,\beta$ small such that $P\in \mathcal{A}_n$ for any $P\in \Sigma$ and $P_{\delta,N}\in\mathcal{A}_{n}^-$.} 
If $H(P_{\delta,N})$ does not tend to zero, we can choose  a subsequence (which we still denote by $P_{\delta,N}$) and  consider $\Sigma_{N}=\Sigma-P_{\delta,N}$. We know from Lemma \ref{lem2.1} that $\Sigma$ has bounded principle curvatures, so the same is true of  $\Sigma_{N}$. Then a subsequence of the $\Sigma_{N}$ will converge smoothly to $\Sigma_\infty$, which is again a  mean convex translating soliton with $H(0)>0$. However, 
\[\inf_{\Sigma_\infty}\frac{\mu^n}{H-\mu^n}=\frac{\mu^n}{H-\mu^n}(0)\leq -\varepsilon_0\]
This contradicts the fact that $Q_\delta$ has no interior negative minimum.

Therefore, we must have $H(P_{\delta,N})\to 0$ as $N\to \infty$. 
\begin{claim}\label{cl2}
If $n\geq 3$, then $\kappa_i(P_{\delta,N})\to 0$, for any $1\leq i\leq n$ as $N\to \infty$. Moreover, there exists $C(n)>1$ such that for all $l$,
\begin{align}
   C(n)^{-1}\leq  \left|\frac{\kappa_l}{H-\mu^n}\right|(P_{\delta,N})\leq C(n)\label{eq:klHmu}\\
    C(n)^{-1}\leq  \frac{H}{H-\mu^n}(P_{\delta,N})\leq C(n).
\end{align}
provided $\delta$ small and $N$ large enough.
\end{claim}
\begin{proof}
Order the principal curvatures as before:
$\kappa_1\leq \cdots\leq \kappa_n$.
 Since $\Sigma$ is uniformly 2-convex,  
\[H=\sum_{j=1}^n\kappa_j\geq \sum_{j=3}^n\kappa_j+\beta H.\]
Therefore $\kappa_j(P_{\delta,N})\to0$ for $j\geq 3 $ and $(\kappa_1+\kappa_2) (P_{\delta,N})$ as $N\to \infty$. By Lemma \ref{lem2.1} ii., we also have $\kappa_j(P_{\delta,N})\to 0,\, j=1,2$ at the same time. Using Lemma \ref{lem2.1} iv., at each $P_{\delta,N}$,
\[0\geq\frac{\mu^n}{H}\geq \frac{\kappa_1-c(n)\sqrt{\delta}\,H }{H}\geq -c({n})\sqrt{\delta} -\frac{1}{n-2}\]
which implies
\begin{align*}
    1\geq \frac{H}{H-\mu^n}=\frac{1}{1-\mu^n/H}\geq \frac{n-2}{n-1+(n-2)c({n})(\delta)}>0
\end{align*}
Suppose $\delta$ is small enough such that at each $P_{\delta,N}$, \eqref{eq:k1aH} holds. Then if $N$ is large enough,
\[-1\leq \frac{\kappa_1}{H-\mu^n}\leq \frac{\mu^n+c(n)\sqrt{\delta}H}{H-\mu^n}\leq Q_{\delta }+c(n)\sqrt{\delta}\leq -\frac12\varepsilon_0+c(n)\sqrt{\delta}\leq -\frac14\varepsilon_0.\]
For $l\geq 2$, it is easy to see
\[\frac{\beta H}{H-\mu^n}\leq \frac{\kappa_l}{H-\mu^n}\leq 1.\]
Our Claim \ref{cl2} follows from taking $\delta$ small enough and $N$ large enough
\end{proof}
% Taking a subsequence of $P_{\delta,N}$ (still denoted as $P_{\delta,N}$) such that $Q_\delta(P_{\delta,N})\to -c_0\in (-1,-\varepsilon_0]$ that is
% \begin{align}
%     \frac{\mu^n}{H-\mu^n}\to -c_0,\quad \frac{H}{H-\mu^n}\to 1+c_0>0,\quad \text{as }N\to \infty
% \end{align}
% then there exists $C>0$ such that, as $N\to \infty$,
It follows from  Claim \ref{cl2} that
\begin{align*}
    -\frac{\mu^n\nabla_l H}{(H-\mu^n)^2}=\frac{-\mu^n}{H-\mu^n}\frac{\kappa_l\langle \tau_l,e_{n+1}\rangle}{H-\mu^n}
\end{align*}
is uniformly bounded. Furthermore, by \eqref{eq:OY} and
\begin{align*}
    \nabla Q_\delta=\frac{(H-\mu^n)\nabla \mu^n-\mu^n(\nabla H-\nabla \mu^n)}{(H-\mu^n)^2}=-\frac{\mu^n\nabla H}{(H-\mu^n)^2}+\frac{H\nabla\mu^n}{(H-\mu^n)^2},
\end{align*}
we see that $\frac{H\nabla\mu^n}{(H-\mu^n)^2}$ is uniformly bounded at each point $P_{\delta,N}$. Since $\frac{H}{H-\mu^n}$ is bounded away from $0$ by Claim \ref{cl2}, we conclude that $\frac{\nabla \mu^n}{H-\mu^n}$ is also uniformly bounded at each point $P_{\delta,N}$. Adopting the notation $\nabla_l=\nabla_{\tau_l}$, 
\begin{claim}\label{cl3} We have
\[\frac{\nabla_{l}H}{H-\mu^n}(P_{\delta,N})\to 0\quad \text{ as } N\to \infty\]
for any $1\leq l\leq n$.
\end{claim}
\begin{proof}
It follows from \eqref{eq:Q-ellip} and \eqref{eq:OY} that
\begin{align}\label{eq:OY-ineq}
    -\frac{2}{N}+2\frac{\left\langle{\nabla (H-\mu^n)}, \nabla Q_\delta\right\rangle}{H-\mu^n}\leq \frac{HF^{ij,pq}h_{ijk}h_{pqk}}{(H-\mu^n)^2}\leq 0
\end{align}
holds at each $P_{\delta,N}$. From the previous analysis, one knows 
\[\frac{\nabla(H-\mu^n)}{H-\mu^n}(P_{\delta,N})\]
is uniformly bounded as $N\to \infty$. Therefore, the left hand side of \eqref{eq:OY-ineq} converges to 0 as $N\to \infty$. Consequently ( since $\frac{H}{H-\mu^n}$ is bounded away from zero)
\begin{align*}
    \frac{F^{ij,pq}h_{ijk}h_{pqk}}{H-\mu^n}\to 0,\quad\text{as } N\to \infty.
\end{align*}
By \eqref{eq:second-deri}, one has
\begin{align}\label{eq:Fto0}
    \frac{F^{ij,pq}h_{ijk}h_{pqk}}{H-\mu^n}\leq \frac{1}{H-\mu^n}\sum_{\{i,j:\kappa_i\neq \kappa_j\}}\sum_{k}\frac{\mu^n_i-\mu^n_j}{\kappa_i-\kappa_j}h_{ijk}^2\leq 0
\end{align}
and each term in the summation is nonpositive. Therefore  all the terms converge to zero as $N\to \infty$.
Suppose $\delta$ is small enough. 
Taking $i=1$ and $j>1$ and using Lemma \ref{lem2.2} ii., one obtains $(H-\mu^n)^{-1}{|h_{1jk}|}\to 0$ as $N\to \infty$ for any $k$. In particular,
\begin{align}\label{eq:h1jjto0}
    \frac{|h_{1jj}|+|h_{11j}|}{H-\mu^n}\to 0,\quad\text{as } N\to \infty,\quad \text{for }j>1.
\end{align}
Taking both $i,j>1$ and $i\geq j$, if $\kappa_i=\kappa_j$, then $h_{ijk}=0$ by \eqref{eq:hijk=0}. If $\kappa_i\neq \kappa_j$, we use Lemma \ref{lem2.3} to conclude
\begin{align*}
    \frac{\Lambda^{n-2}}{2\cdot n!}\frac{h_{ijk}^2}{(H-\mu^n)^2}\leq& \frac{\Lambda^{n-2}}{2\cdot n!}\frac{1}{H-\mu^n}\frac{h_{ijk}^2}{|\kappa_i-\kappa_j|+\sqrt{\kappa_i\kappa_j}} \\
    \leq&-\frac{1}{H-\mu^n}\frac{\mu^i-\mu^j}{\kappa_i-\kappa_j}h_{ijk}^2\leq \left|\frac{F^{ij,pq}h_{ijk}h_{pqk}}{H-\mu^n}\right|
\end{align*}
In either case, we must have 
\begin{align}\label{eq3.20}
    \frac{|h_{ijk}|}{H-\mu^n}\to 0\quad \text{as }N\to \infty, \quad\text{for }i\neq j.
\end{align}
Since $\nabla_l \mu^n=\sum_i\mu^n_i\nabla_l\kappa_i=\sum_i\mu^n_ih_{iil}$, we can rewrite $\nabla_l Q_{\delta}$ as 
\begin{eqnarray}\label{eq3.30}
    \nabla_l Q_\delta=\frac{H\nabla_l\mu^n-\mu^n\nabla_lH}{(H-\mu^n)^2}&=\left[1+\frac{H(\mu^n_l-1)}{H-\mu^n}\right]\frac{h_{lll}}{H-\mu^n}\\
\nonumber   +& \sum_{i\neq l}\frac{H\mu^n_lh_{iil}-\mu^n h_{iil}}{(H-\mu^n)^2}
\end{eqnarray}
When $\delta$ is small enough,
it follows from Lemma \ref{lem2.2} ii.  that $\mu_l^n\to 0 $  if $l>1$ and $\mu_l^n\to  1$ if $l=1$, as $N\to \infty$, uniformly as $\delta$ decreases to zero. Combined with the uniform estimates of Claim \ref{cl2} and \eqref{eq3.20}, the last term on the right hand side of \eqref{eq3.30} tends to zero as $N\to \infty$. Moreover,
  $1+H(\mu^n_l-1)(H-\mu^n)^{-1}$ is uniformly bounded away from zero for any $l$. Since $\nabla_lQ_\delta\to 0$ as $N\to \infty$, we conclude  from \eqref{eq3.30} that 
  \be \label{eq3.40}
   \frac{h_{lll}}{H-\mu^n} \to 0\,\,\text{as}\,\, N\to \infty.
   \ee
Therefore \eqref{eq3.20} and \eqref{eq3.40} imply
\[\frac{\nabla_lH}{H-\mu^n}=\sum_{i=1}^n\frac{h_{iil}}{H-\mu^n}\to 0
\,\,\text{as}\,\, N\to \infty,\]
proving Claim \ref{cl3}.
% Consider $(H-\mu)^{-1}\nabla_1\mu$. Since $\nabla_1\mu=\sum_{i=1}^n\mu^i h_{ii1}$, then
% It follows from \eqref{eq:h1jjto0} and \eqref{eq:OY} that $(H-\mu)^{-1}|h_{111}|\to 0$ as $N\to \infty$. Hence, it is easy to verify that $(H-\mu)^{-1}\nabla_1\mu\to 0$ as $N\to \infty$.

% Now let us consider $(H-\mu)^{-1}\nabla_l\mu$ for $l>1$. Let us just take $l=2$ for example, the other cases are similar. 
% \begin{align*}
%     \frac{\nabla_2\mu}{H-\mu}=\frac{\mu^1h_{112}}{H-\mu}+\frac{\mu^2h_{222}}{H-\mu}+\sum_{j=3}^n\frac{\mu^jh_{jj2}}{H-\mu}
% \end{align*}

\end{proof}

Again consider $\Sigma_N=\Sigma-P_{\delta,N}$. A subsequence of the $\Sigma_N$ converges locally smoothly to a translator $\Sigma_\infty$ .
Since $H=\langle\nu, e_{n+1}\rangle$, 
\[\frac{\nabla_lH}{H-\mu^n}(P_{\delta,N})=\frac{\kappa_l\langle\tau_l,e_{n+1}\rangle}{H-\mu^n}(P_{\delta,N})\to 0,\quad  \text{ as }N\to \infty.\]
But then by \eqref{eq:klHmu},
\begin{align}
     \langle\tau_l,e_{n+1}\rangle(P_{\delta,N})\to 0\text{ as }N\to \infty.
\end{align}
This means that at the origin on $\Sigma_{\infty}$, we have $\langle\tau_l,e_{n+1}\rangle(0)=0$. In other words, $\nu=e_{n+1}$ and $H=1$ at $0$.
However, 
\[\inf_{\Sigma_\infty}\frac{\mu^n}{H-\mu^n}=\frac{\mu^n}{H-\mu^n}(0)\leq -\varepsilon_0\]
This again contradicts the fact that $Q_\delta$ cannot have an interior negative minimum. This completes the proof of Theorem \ref{thm3.1} and consequently Theorem \ref{thm1.1} is proven.

\end{proof}

\bibliographystyle{plainnat}
\bibliography{2-convex}

\end{document}